\documentclass[12pt,]{amsart}
\begin{document}
\bibliographystyle{tom}

\newtheorem{lemma}{Lemma}
\newtheorem{thm}[lemma]{Theorem}
\newtheorem{cor}[lemma]{Corollary}
\newtheorem{voorb}[lemma]{Example}
\newtheorem{rem}[lemma]{Remark}
\newtheorem{prop}[lemma]{Proposition}  
\newtheorem{deff}[lemma]{Definition}  

\newtheorem{stat}[lemma]{{\hspace{-5pt}}}

\newenvironment{remarkn}{\begin{rem} \rm}{\end{rem}}
\newenvironment{exam}{\begin{voorb} \rm}{\end{voorb}}

\newcommand{\esssup}{\mathop{\rm ess\,sup}}
\newcommand{\essinf}{\mathop{\rm ess\,inf}}
\newcommand{\divv}{\mathop{\rm div}}
\newcommand{\supp}{\mathop{\rm supp}}

\newcommand{\cl}{{\cal L}}

\makeatletter
\def\eqnarray{\stepcounter{equation}\let\@currentlabel=\theequation
\global\@eqnswtrue
\tabskip\@centering\let\\=\@eqncr
$$\halign to \displaywidth\bgroup\hfil\global\@eqcnt\z@
  $\displaystyle\tabskip\z@{##}$&\global\@eqcnt\@ne
  \hfil$\displaystyle{{}##{}}$\hfil
  &\global\@eqcnt\tw@ $\displaystyle{##}$\hfil
  \tabskip\@centering&\llap{##}\tabskip\z@\cr}

\def\endeqnarray{\@@eqncr\egroup
      \global\advance\c@equation\m@ne$$\global\@ignoretrue}
\makeatother

\newcounter{teller}
\renewcommand{\theteller}{\Roman{teller}}
\newenvironment{tabel}{\begin{list}%
{\rm \bf \Roman{teller}.\hfill}{\usecounter{teller} \leftmargin=1.1cm
\labelwidth=1.1cm \labelsep=0cm \parsep=0cm}
                      }{\end{list}}

\newcommand{\Ri}{{\bf R}}
\newcommand{\ch}{{\cal H}}
\newcommand{\cd}{{\cal D}}

\newcommand{\one}{1\hspace{-4.5pt}1}

\hyphenation{Lipschitz}

\title{ Ellipticity and Ergodicity} 

\author{Derek W. Robinson and Adam Sikora}

\address{Derek W. Robinson, 
Centre for Mathematics and its Applications,
Mathematical Sciences Institute,
Australian National University,
Canberra, ACT 0200,
Australia}
\email{Derek.Robinson@anu.edu.au}
\address{Adam Sikora, Department of Mathematical Sciences, New Mexico State
University, Las Cruces, NM 88003-8001, USA
{\rm and} Department of Mathematics,
Australian National University,
ACT 0200 Australia}
\email{asikora@nmsu.edu\textrm{, } sikora@maths.anu.edu.au}

\begin{abstract}
Let $S=\{S_t\}_{t\geq0}$ be the submarkovian semigroup  on $L_2(\Ri^d)$ generated by a self-adjoint, second-order, 
divergence-form, elliptic operator $H$  with Lipschitz continuous coefficients $c_{ij}$.
Further let $\Omega$ be an open subset of $\Ri^d$.
Under the  assumption that $C_c^\infty(\Ri^d)$ is a core for $H$   we prove that  $S$ leaves $L_2(\Omega)$ invariant if, and only if,
it is invariant under the flows generated by the vector fields $Y_i=\sum^d_{j=1}c_{ij}\partial_j$.
\end{abstract}

\date{December 2007}

\subjclass{35J70, 35Hxx, 35F05, 31C15}

\maketitle

 \section{Introduction}\label{S1}

Let $S$ be a submarkovian semigroup on $L_2(\Ri^d)$ generated by a self-adjoint second-order elliptic operator $H$ in divergence form. 
If the operator is strongly elliptic  then $S$ acts ergodically, i.e.\ there are no non-trivial $S$-invariant subspaces of $L_2(\Ri^d)$.
Nevertheless there are many examples of degenerate elliptic operators for which there are 
subspaces $L_2(\Omega)$  invariant under the action of $S$ (see, for example, \cite{ERSZ2}
\cite{RSi} \cite{ER29}).
Our aim is to examine operators with coefficients  which are  Lipschitz continuous and  characterize the  $S$-invariance of $L_2(\Omega)$ by the invariance under a family of associated  flows. 
Then one can combine the characterization with a domination estimate to establish invariance properties for a large class of degenerate operators with $L_\infty$-coefficients.
In order to formulate our main result more precisely we need  some further notation.

First  define $H$ as the Friederichs extension of the positive symmetric operator $H_0$ with
 domain $D(H_0)=C_c^\infty(\Ri^d)$ and action
\[
H_0\varphi= -\sum^d_{i,j=1}\partial_i\,c_{ij}\partial_j\varphi
\]
where the coefficients $c_{ij}=c_{ji}\in W_{\rm loc}^{1,\infty}(\Ri^d)$ are real and  $C=(c_{ij})$
is a  positive-definite matrix over $\Ri^d$.
Thus $H$ is a divergence form  operator with coefficients which are locally Lipschitz continuous.
It is the positive, self-adjoint, operator on $L_2(\Ri^d)$ associated with the closure $\overline h_0$  of the quadratic form 
\begin{equation}
h_0(\varphi)= \sum^d_{i,j=1}(\partial_i\varphi, c_{ij}\partial_j\varphi)
\label{esep1.1}
\end{equation}
with domain $ D(h_0)=C_c^\infty(\Ri^d)$.
Set $h=\overline h_0$.
Then  $h$ is  a Dirichlet form and the self-adjoint semigroup $S$ generated by $H$ is automatically submarkovian
(for details on Dirichlet forms and submarkovian semigroups see \cite{FOT} or \cite{BH}).

Secondly, let $\psi\in C_c^\infty(\Ri^d)$ and define $Y_\psi$ as the $L_2$-closure of the first-order
partial differential operator  with action
\[
Y_\psi\varphi= \sum^d_{i,j=1}c_{ij}(\partial_i\psi)\partial_j\varphi
\]
for all $\varphi\in C_c^\infty(\Ri^d)$.
The   $Y_\psi$ generate positive, continuous,  one-parameter quasi-contractive groups $ T^\psi$ on $L_2(\Ri^d)$.  
The latter result follows because the coefficients $c_j=\sum^d_{i=1}c_{ij}(\partial_i\psi)$ of the  $Y_\psi$ are in $W^{1,\infty}(\Ri^d)$ (see Lemma~\ref{gllc} for details).
We refer to these groups  as flows.

Our primary aim is to establish the following characterization of invariance.

\begin{thm}\label{tsep2} 
Assume  $c_{ij}\in W_{\rm loc}^{1,\infty}(\Ri^d)$ and that $C_c^\infty(\Ri^d)$ is a core for $H$.
Let $\Omega$ be a measurable subset of $\Ri^d$.

The following conditions are equivalent.
\begin{tabel}
\item\label{tsep2-1}
$S_tL_2(\Omega)\subseteq L_2(\Omega)$ for all $t>0$.
\item\label{tsep2-2}
$T^\psi_tL_2(\Omega)=L_2(\Omega)$ for all $\psi\in C_c^\infty(\Ri^d)$ and all $t\in\Ri$. 
\end{tabel}
\end{thm}

If the coefficients $c_{ij}\in W^{1,\infty}(\Ri^d)$  one can characterize the $S$-invariance of $L_2(\Omega)$ with the  flows generated by the $L_2$-closures of the operators $Y_i$  with  action
\[
Y_i\varphi= \sum^d_{j=1}c_{ij}\partial_j\varphi
\]
for all $\varphi\in C_c^\infty(\Ri^d)$ where $i\in\{1,\ldots,d\}$.

\begin{thm}\label{tsep1} 
Assume  $c_{ij}\in W^{1,\infty}(\Ri^d)$ and that $C_c^\infty(\Ri^d)$ is a core for $H$.
Let $\Omega$ be a measurable subset of $\Ri^d$.

The following conditions are equivalent:
\begin{tabel}
\item\label{tsep1-1}
$S_tL_2(\Omega)\subseteq L_2(\Omega)$ for all $t>0$,
\item\label{tsep1-2}
$T^{(i)}_tL_2(\Omega)=L_2(\Omega)$ for all  $i\in\{1,\ldots,d\}$ and all $t\in\Ri$.
\end{tabel}
\end{thm}

The operators $Y_i$  were used by  Oleinik and Radkevic \cite{OR} to analyze 
 hypoellipticity and subellipticity properties
of degenerate elliptic operators $H$ with $C^\infty$-coefficients $c_{ij}$.
In particular they established that if the $Y_i$ satisfy H\"ormander's condition \cite{Hor1}, i.e.\ if the $Y_i$ and their multi-commutators up to some finite order $k$ span the tangent space at each point $x\in\Ri^d$, then $H$ satisfies the subelliptic estimate  $H\geq \mu\,\Delta^\gamma-\nu\,I$ for some $\mu>0$, $\nu\geq0$ and $\gamma\in\langle0,1]$ where $\Delta$ is the usual Laplacian.
There is, however, no simple relationship between the values of $k$ and $\gamma$
(see \cite{JSC} for a review of related results).

The foregoing characterizations of $S$-invariant subspaces differ from the earlier results \cite{RSi} \cite{ER29} based on capacity estimates on the boundary $\partial\Omega$ of $\Omega$. Here our results  require 
the assumption that $C_c^\infty(\Ri^d)$ is a core for $H$. Unfortunately, this assumption does not always holds for 
degenerate elliptic operators, see \cite{RSi09}. We do not know if this supposition  is indeed necessary.

\section{Preliminaries}\label{S2}

The proofs  of Theorems~\ref{tsep2} and \ref{tsep1}  depend on some basic properties of  flows
 which follow from standard results on first-order partial differential equations which we first summarize.

Let  $c_j\in W^{1,\infty}(\Ri^d)$ for each  $j\in\{1,\ldots,d\}$ and let 
$Y=\sum^d_{j=1}c_j\partial_j$ be the corresponding first-order operator on $C_1(\Ri^d)$.
The operator  $Y$ is a model for the operators $Y_\psi$  and $Y_i$ introduced in Section~\ref{S1}.
Since  $c_j\in W^{1,\infty}(\Ri^d)$
  there exists a unique Lipschitz continuous function $(t,x)\in \Ri \times \Ri^d\mapsto f_t(x)\in\Ri^d$
satisfying the evolution equation
\begin{equation}
\frac{\partial f_t(x) }{\partial t} = c( f_t(x))
\;\;\;\; \mbox{and} \;\;\;\; 
f_0(x) = x
\label{elc3.0}
\end{equation}
for all  $x \in \Ri^d$ and  $t \in \Ri$ (see, for example, \cite{Hille} Chapter~2 and in particular Theorem~2.6.3).
Moreover, there are $M,\omega\geq0$ such that 
\[
|f_t(x)-f_s(x)|\leq M\,|t-s|\,e^{\omega(|t|\wedge|s|)}
\]
and 
\[
|f_t(x)-f_t(y)|\leq |x-y|\,e^{\omega|t|}
\]
for all $x,y\in\Ri^d$ and $s,t\in\Ri$ (see \cite{Hille}, Theorems~3.1.1 and 3.2.1).
The partial derivatives $\partial_i(f_t)_j$ of the components $(f_t)_j$ of $f_t$ 
and  the Jacobian $J_t$ of the transformation $x\in\Ri^d\to f_t(x)$ are bounded.
Moreover,  there are 
$M'\geq 1$ and $\omega'\geq0$ such that 
\begin{equation}
\|J_t\|_\infty=\|\det(\partial_i(f_t)_j)\|_\infty\leq M'\,e^{\omega'|t|}
\label{elc3.21}
\end{equation}
for all $t\in\Ri$.
We   adopt the conventional notation $\exp(tY)(x) = f_t(x)$
and then one has the group property \[
(\exp(tY)(\exp(sY))(x) = \exp((t+s)Y)(x)
\]
 for all $x \in \Ri^d$ and $t,s \in \Ri$.
 
 One can immediately define a  positive, continuous,  one-parameter group 
 of isometries $T_t $ on $ C_b(\Ri^d)$  by setting 
$(T_t \psi)(x) = \psi(\exp(tY)(x))=(\psi\circ f_t)(x)$.
In particular $T$ is conservative, i.e.\  $T_t\one =\one$ for all $t\in\Ri$.
Moreover, the group $T$ extends to a conservative weakly$^*$ continuous group of isometries on $L_\infty(\Ri^d)$.
But the Jacobian of the transformation $x\to e^{tY}x$ is uniformly bounded  by (\ref{elc3.21}).
Therefore $T$ also extends to a strongly continuous one-parameter group on the spaces $L_p(\Ri^d)$
for $p\in[1,\infty\rangle$.
In fact one has the following statement  in which  the weak$^*$ topology is to be understood if $p=\infty$.

\begin{lemma}\label{gllc}
The $L_p$-closure  of $Y|_{C_c^\infty(\Ri^d)}$ generates a positive,  continuous,  one-parameter group $T^{(p)}$ on $L_p(\Ri^d)$, for each $p\in[1,\infty]$, with action
\begin{equation}
(T^{(p)}_t \psi)(x) = \psi(\exp(tY)(x))
\label{elc3.30}
\end{equation}
for  $\psi\in L_p(\Ri^d)$ and $t\in\Ri$.
Moreover, $T^{(p)}_tW^{1,p}(\Ri^d)=W^{1,p}(\Ri^d)$ and
\begin{equation}
\|T^{(p)}_t\|_{p\to p}\leq e^{\nu|t|/p}
\label{elc3.31}
\end{equation}
for all $t\in\Ri$ with $\nu=\|\divv c\,\|_\infty$.
\end{lemma}
\proof\  First it follows from the definition (\ref{elc3.30}) that $T^{(p)}$ is a positive, strongly continuous, group
on $L_p(\Ri^d)$ if $p\in[1,\infty\rangle$ and $T^{(\infty)}$ is a positive, weakly$^*$ continuous, group of 
isometries on $L_\infty(\Ri^d)$.

Secondly, $W^{1,p}(\Ri^d)$  is $T^{(p)}$-invariant because the bounds on the partial derivatives
$\partial_i(f_t)_j$ immediately give bounds
\[
\|\partial_iT^{(p)}_t\psi\|_p\leq M\, e^{\omega|t|}\,\sup_{1\leq j\leq d}\|\partial_j\psi\|_p
\]
for all $\psi\in W^{1,p}(\Ri^d)$.

Thirdly, it follows from the definition (\ref{elc3.30}) that the generator of $T^{(p)}$ is a closed extension
of the operator $Y$ defined on $W^{1,p}(\Ri^d)$.
But since $W^{1,p}(\Ri^d)$  is $T^{(p)}$-invariant it must be a core of the generator.
Hence the generator $Y_p$  is the $L_p$-closure  of $Y|_{W^{1,p}(\Ri^d)}$.
Then   $Y_p$ is the $L_p$-closure  of $Y|_{C_c^\infty(\Ri^d)}$
since $C_c^\infty(\Ri^d)$ is dense in $W^{1,p}(\Ri^d)$.

Fourthly, the  adjoint  $Y^*_\infty$ of $Y_\infty$ generates a strongly continuous group of isometries on $L_1(\Ri^d)$.
But $Y^*_\infty\psi=-Y\psi-(\divv c)\psi$ for all $\psi\in W^{1,1}(\Ri^d)$. 
Therefore   $Y_\infty^*\supseteq -Y_1-(\divv c)I$ and since both operators are generators of continuous groups one must have   $Y_\infty^*= -Y_1-(\divv c)I$.
Then the group $T^{(1)}$ generated by $Y_1=-Y^*_\infty-(\divv c)I$ 
  satisfies the bounds  $\|T^{(1)}_t\|_{1\to1}\leq \exp(\nu|t|)$
for all $t\in\Ri$  by the Trotter product formula.
The $L_p$-bounds (\ref{elc3.31})  follow by interpolation.
\hfill$\Box$

The functions $L_\infty(\Ri^d)$ act as multipliers on the $L_p$-spaces and it follows from the
 action~(\ref{elc3.30}) of the group that
\begin{equation}
T^{(p)}_t(\varphi\psi)=(T^{(\infty)}_t \varphi)(T^{(p)}_t\psi)
\label{eeqn}
\end{equation}
for all $\varphi\in L_\infty(\Ri^d)$ and $\psi\in L_p(\Ri^d)$.

\begin{lemma}\label{linv1}
Let $\Omega$ be a measurable subset of $ \Ri^d$.
The following conditions are equivalent.
\begin{tabel}
\item\label{linv1-1}
$T^{(p)}_tL_p(\Omega)=L_p(\Omega)$ for all $t\in\Ri$ and for one $($for all$)$ $p\in[1,\infty\rangle$.
\item\label{linv1-10}
$[T^{(p)}_t, \one_\Omega]=0$ for all $t\in\Ri$ and for one $($for all$)$ $p\in[1,\infty\rangle$.
\item\label{linv1-2}
$T^{(\infty)}_t\one_\Omega=\one_\Omega$ for all $t\in\Ri$.
\end{tabel}
\end{lemma}
\proof\
First suppose Condition~\ref{linv1-1} is valid for one $p\in[1,\infty\rangle$.
Then if $q\in[1,\infty\rangle$ one has
\[
T^{(q)}_t(L_p(\Omega)\cap L_q(\Omega))=T^{(p)}_t(L_p(\Omega)\cap L_q(\Omega))
\subseteq L_p(\Omega)\cap L_q(\Ri^d)=L_p(\Omega)\cap L_q(\Omega)
\]
for all $t\in\Ri$.
Then since $L_p(\Omega)\cap L_q(\Omega)$ is dense in $L_q(\Omega)$ one has 
$T^{(q)}_t L_q(\Omega)\subseteq  L_q(\Omega)$ for all $t\in\Ri$ and by the group property
$T^{(q)}_t L_q(\Omega)= L_q(\Omega)$ for all $t\in\Ri$.
Thus Condition~\ref{linv1-1} is valid for all $p\in[1,\infty\rangle$.

Similarly, since $[T^{(p)}_t,\one_\Omega]\psi=[T^{(q)}_t,\one_\Omega]\psi$ for all $\psi\in L_p(\Ri^d)\cap L_q(\Ri^d)$ it follows that if Condition~\ref{linv1-10} is valid for one  $p\in[1,\infty\rangle$
then it is valid for all $p\in[1,\infty\rangle$.

\smallskip

Now it suffices to prove the equivalence of the three conditions with $p=2$.

\smallskip

\noindent\ref{linv1-1}$\Rightarrow$\ref{linv1-10}
Since $\one_\Omega$ is the orthogonal projection from $L_2(\Ri^d)$ onto the invariant subspace $L_2(\Omega)$ it follows that
\[
T^{(2)}_t\one_\Omega=\one_\Omega\, T^{(2)}_t\one_\Omega
\;.
\]
Then by taking adjoints 
\[
\one_\Omega \,T^{(2)\,*}_t=\one_\Omega \,T^{(2)\,*}_t\one_\Omega
\;.
\]
But since the generator $Y_2$ of $T^{(2)}$ satisfies $Y_2=-Y_2^*-(\divv c)I$ it follows by another application
of the Trotter product formula that 
\[
\one_\Omega \,T^{(2)}_t=\one_\Omega \,T^{(2)}_t\one_\Omega
\;.
\]
Therefore $[T^{(2)}_t,\one_\Omega]=0$.

\smallskip

\ref{linv1-10}$\Rightarrow$\ref{linv1-2}
It follows immediately from Condition~\ref{linv1-10} and (\ref{eeqn}) that 
\[
\one_\Omega \,T^{(2)}_t\psi= T^{(2)}_t\one_\Omega\psi
=(T^{(\infty)}_t\one_\Omega)T^{(2)}_t\psi
\]
for all $\psi\in L_2(\Ri^d)$.
Replacing $\psi$ by $T^{(2)}_{-t}\psi$ one deduces that $T^{(\infty)}_t\one_\Omega=\one_\Omega$.

\smallskip

\ref{linv1-2}$\Rightarrow$\ref{linv1-1}
If $\psi\in L_2(\Ri^d)$ then
\[
T^{(2)}_t\one_\Omega\psi
=(T^{(\infty)}_t\one_\Omega)T^{(2)}_t\psi=\one_\Omega \,T^{(2)}_t\psi
\]
by (\ref{eeqn}) and Condition~\ref{linv1-2}.
Therefore $T^{(2)}_tL_2(\Omega)\subseteq L_2(\Omega)$ for all $t\in\Ri$.
Since $T^{(2)}$ is a group  it then follows that 
$T^{(2)}_tL_2(\Omega)= L_2(\Omega)$ for all $t\in\Ri$.
\hfill$\Box$

\bigskip

Lemma~\ref{linv1} has the following straightforward corollary.

\begin{cor}\label{cinv}
The following conditions are equivalent.
\begin{tabel}
\item\label{cinv1-1}
$T^{(2)}_tL_2(\Omega)=L_2(\Omega)$ for all $t\in\Ri$.
\item\label{cinv1-2}
$(Y^*\varphi,\one_\Omega)=0$ for all $\varphi\in C_c^\infty(\Ri^d)$.
\end{tabel}
\end{cor}
\proof
\ref{cinv1-1}$\Rightarrow$\ref{cinv1-2}
It follows from Lemma~\ref{linv1} that Condition~\ref{cinv1-1} is equivalent to
$T^{(\infty)}_t\one_\Omega=\one_\Omega$ for all $t\in\Ri$ and this is clearly equivalent to
\[
(T^{(\infty)\,*}_t\varphi, \one_\Omega)=(\varphi,\one_\Omega)
\]
for all $t\in\Ri$ and $\varphi\in L_1(\Ri^d)$.
Then Condition~\ref{cinv1-2} follows by differentiation.

\smallskip

\noindent\ref{cinv1-2}$\Rightarrow$\ref{cinv1-1}
Since $Y^*=-Y-(\divv c)I$ and $C_c^\infty(\Ri^d)$ is a core of $Y_1$ it follows 
from Condition~\ref{cinv1-2}  by closure that
$((Y_1+(\divv c)I)\varphi, \one_\Omega)=0$ for all $\varphi\in D(Y_1)$.
But $(Y_1+(\divv c)I)=-Y^*_\infty$.
Therefore Condition~\ref{cinv1-2} implies 
$(Y^*_\infty\varphi,\one_\Omega)=0$ for all $\varphi\in D(Y^*_\infty)=D(Y_1)$.
Now Duhamel's formula gives
\[
(T^{(\infty)\,*}_t\varphi, \one_\Omega)-(\varphi,\one_\Omega)
=\int^t_0ds\,(Y^*_\infty \,T^{(\infty)\,*}_s\varphi, \one_\Omega)=0
\]
for all $\varphi\in D(Y^*_\infty)$.
Then since $D(Y^*_\infty)=D(Y_1)$ is dense in $L_1(\Ri^d)$ it follows that  $T^{(\infty)}_t\one_\Omega=\one_\Omega$ for all $t\in\Ri$.
Condition~\ref{cinv1-1}  follows from Lemma~\ref{linv1}.
\hfill$\Box$

\section{Proofs of Theorems~\ref{tsep2} and \ref{tsep1} }\label{S3}

First recall that $h$ denotes the Dirichlet form associated with the elliptic operator $H$, i.e.\ $h$ is the closure of the form $h_0$ defined by  (\ref{esep1.1}).
The form $H$ is local in the sense  of \cite{FOT}, i.e.\  if $\varphi,\psi\in D( h)$  and $\varphi\psi=0$  then  the sesquilinear form associated with $h$ satisfies $h(\varphi,\psi)=0$.
(This condition appears somewhat  stronger than that of \cite{FOT} but it is 
 in fact equivalent 
by a result of  Schmuland \cite{Schm}.)
Secondly,  it follows from the Dirichlet property that $D(h)\cap L_\infty(\Ri^d)$ is an algebra.
Therefore,  for each positive $\xi\in D( h)\cap L_\infty(\Ri^d)$,  one can define the truncation $h_\xi$ of $h$ by
 \begin{equation}
 h_{\xi}(\varphi,\psi)=2^{-1}\left(h(\xi\varphi,\psi)+ h(\varphi,\xi\psi)
- h(\xi,\varphi\psi)\right)
\label{ep}
\end{equation}
for all $\varphi,\psi\in  D( h)\cap L_\infty(\Ri^d)$.
The Dirichlet property  implies that  $D(h)\subseteq D(h_\xi)$ and 
\[
0\leq  h_\xi(\varphi)\leq \|\xi\|_\infty\,h(\varphi)
\]
 for all $\varphi\in D( h)$ (see, for example, \cite{BH}, Proposition~4.1.1).
 The truncated form $ h_\xi$ is not necessarily closed but it is local. 
The locality of  $ h_{\xi}$ follows straightforwardly from the locality of $ h$.

Thirdly,  for each $\psi\in C_c^\infty(\Ri^d)$ let $Y_\psi$ denote the vector field on $\Ri^d$ defined in Section~\ref{S1}
  The operators $Y_\psi$ are related to the truncated forms $h_\xi$ by the identity
 \begin{equation}
 (\xi,Y_\psi\varphi)=\int_{\Ri^d}dx\,\xi(x)\sum^d_{i,j=1}c_{ij}(x)(\partial_j\psi)(x)(\partial_i\varphi)(x)
 =h_\xi(\psi,\varphi)
 \label{elc3.1}
 \end{equation}
which is certainly valid for all for all $\varphi,\psi\in C_c^\infty(\Ri^d)$ and positive $\xi\in D( h)\cap L_\infty(\Ri^d)$.
But
\begin{eqnarray*}
| (\xi,Y_\psi\varphi)|&\leq &\|\xi\|_\infty\,\|Y_\psi\varphi\|_1\\[5pt]
&=&\|\xi\|_\infty \int_{\Ri^d}dx\,\bigg|\sum^d_{i,j=1}c_{ij}(x)(\partial_i\psi)(x)(\partial_j\varphi)(x)\bigg|\\[5pt]
&\leq&\|\xi\|_\infty \int_{\Ri^d}dx\,\bigg(\sum^d_{i,j=1}c_{ij}(x)(\partial_i\psi)(x)(\partial_j\psi)(x)\bigg)^{1/2}\,
\bigg(\sum^d_{i,j=1}c_{ij}(x)(\partial_i\varphi)(x)(\partial_j\varphi)(x)\bigg)^{1/2}\\[5pt]
&\leq&\|\xi\|_\infty\, h(\psi)^{1/2}\,h(\varphi)^{1/2}\;.
\end{eqnarray*}
Therefore (\ref{elc3.1}) extends by continuity to all $\varphi,\psi\in D( h)\cap L_\infty(\Ri^d)$
always with $\xi\in D( h)\cap L_\infty(\Ri^d)$ positive.
Moreover, it follows from this estimation that 
\[
\|Y_\psi\varphi\|_1\leq h(\psi)^{1/2}\, h(\varphi)^{1/2}
\]
for  all $\varphi,\psi\in C_c^\infty(\Ri^d)$  and then by continuity  for all $\varphi\in D(h)$.
In particular $Y_\psi D( h)\subseteq L_1(\Ri^d)$ for all $\psi\in C_c^\infty(\Ri^d)$.

Now we are prepared to prove the first theorem.

\smallskip

\noindent{\bf Proof of Theorem~\ref{tsep2}}
\ref{tsep2-1}$\Rightarrow$\ref{tsep2-2}
Note that the proof of this implication does not use the assumption that $C_c^\infty(\Ri^d)$ is a core of $H$.

It follows from Condition~\ref{tsep2-1} that $\one_\Omega D(h)\subseteq D(h)$ (see \cite{FOT}, Theorem~1.6.1).
Thus if $\xi\in D( h)\cap L_\infty(\Ri^d)$  is positive then $\one_\Omega\xi\in D( h)\cap L_\infty(\Ri^d)$ 
and  $\one_\Omega\xi$ is also positive.
Now it follows from (\ref{elc3.1}) that 
\[
 (\xi,\one_\Omega Y_\psi\varphi) =  (\one_\Omega\xi,Y_\psi\varphi) = h_{\one_\Omega\xi}(\psi, \varphi)
  \]
  for all $\psi\in C_c^\infty(\Ri^d)$ and $\varphi\in D( h)\cap L_\infty(\Ri^d)$.
Then by (\ref{ep}) and locality  
  \begin{eqnarray*}
 h_{\one_\Omega\xi}(\psi, \varphi)&=&
2^{-1}\left(h((\one_\Omega\xi)\psi,\varphi)+ h(\psi,(\one_\Omega\xi)\varphi))
-  h(\one_\Omega\xi,\psi\varphi)\right)\\[5pt]
&=&
2^{-1}\Big(h(\one_\Omega(\xi\psi),\one_\Omega\varphi)+ h(\one_\Omega\psi,\one_\Omega(\xi\varphi))
- h(\one_\Omega\xi,\one_\Omega\psi\varphi)\Big)\\[5pt]
&=&
2^{-1}\Big( h(\xi(\one_\Omega\psi),\one_\Omega\varphi)+ h(\one_\Omega\psi,\xi(\one_\Omega\varphi))
-  h(\xi,(\one_\Omega\psi)(\one_\Omega\varphi))\Big)\\[5pt]
&=&
 h_\xi(\one_\Omega\psi,\one_\Omega\varphi)=h_\xi(\psi,\one_\Omega\varphi)
\;.
\end{eqnarray*}
Now another application of (\ref{elc3.1}) gives
\[
 (\xi,\one_\Omega Y_\psi\varphi)=h_{\one_\Omega\xi}(\psi, \varphi)
 =h_\xi(\psi,\one_\Omega\varphi)=(\xi,Y_\psi\one_\Omega \varphi)
 \]
 for all $\xi, \varphi \in  D( h)\cap L_\infty(\Ri^d)$.
Hence 
\[
\one_\Omega Y_\psi\varphi=Y_\psi\one_\Omega \varphi
\]
for all $\varphi \in D( h)\cap L_\infty(\Ri^d)$ and all $\psi\in C_c^\infty(\Ri^d)$.
In particular this is valid for all $\varphi\in W^{1,2}(\Ri^d)\cap L_\infty(\Ri^d)$.
But it follows from Lemma~\ref{gllc} that $T^\psi_tW^{1,2}(\Ri^d)\cap L_\infty(\Ri^d)=W^{1,2}(\Ri^d)\cap L_\infty(\Ri^d)$
for all $t\in \Ri$.
Therefore the Duhamel formula gives
\[
[T^\psi_t,\one_\Omega]\varphi=-\int^t_0ds\, T^\psi_{t-s}[Y_\psi,\one_\Omega] T^\psi_s \varphi=0
\]
for all $\varphi$ in the $L_2$-dense subspace $W^{1,2}(\Ri^d)\cap L_\infty(\Ri^d)$.
Therefore Condition~\ref{tsep2-2} now follows from Lemma~\ref{linv1}.

\smallskip

\noindent\ref{tsep2-2}$\Rightarrow$\ref{tsep2-1}
 Let $y_\psi$ denote the linear functional
  \[
\varphi\in C_c^\infty(\Ri^d)\mapsto   y_\psi(\varphi)  =\int_{\Ri^d}dx\,(Y_\psi\varphi)(x)
\]
where $\psi\in C_c^\infty(\Ri^d)\subseteq D(H)$ is  held fixed.
  It follows from the definition of $Y_\psi$ that 
  \[
y_\psi(\varphi) =h(\psi, \varphi )=(H\psi,\varphi)
  \]
  for all $\varphi\in C_c^\infty(\Ri^d)$.
  Then since $|(H\psi,\varphi)|\leq \|H\psi\|_\infty\|\varphi\|_1$ one deduces that $y_\psi$ extends
  to a continuous linear  functional over $L_1(\Ri^d)$ satisfying the bounds $| y_\psi(\varphi)|\leq \|H\psi\|_\infty\|\varphi\|_1$.
  
 Condition~\ref{tsep2-2} implies, however, that $T^\psi_t\one_\Omega=\one_\Omega$ on $L_\infty(\Ri^d)$ for all $t\in\Ri$ by Lemma~\ref {linv1}.
 Therefore $T^\psi_t\one_\Omega\varphi=\one_\Omega T^\psi_t\varphi$  for all $t\in\Ri$ on $L_1(\Ri^d)$.
 Thus if $\varphi\in D(Y_\psi)$ then $\one_\Omega\varphi\in D(Y_\psi)$ and 
 $Y_\psi\one_\Omega\varphi=\one_\Omega Y_\psi\varphi$ on $L_1(\Ri^d)$.
 Therefore
 \begin{eqnarray}
(H\psi,\one_\Omega\varphi)
&=&\int_{\Ri^d} dx\,(Y_\psi\one_\Omega\varphi)(x)\nonumber\\[5pt]
&=&\int_{\Ri^d} dx\,(\one_\Omega Y_\psi\varphi)(x)=\int_\Omega dx\,(Y_\psi\varphi)(x)
\label{elc1.1}
\end{eqnarray}
for all $\varphi\in C_c^\infty(\Ri^d)$.
In particular
\[
|(H^{1/2}(H^{1/2}\psi),\one_\Omega\varphi)|=
|(H\psi,\one_\Omega\varphi)|\leq \|Y_\psi\varphi\|_1\leq h(\psi)^{1/2}\, h(\varphi)^{1/2}
\;.
\]
Then since $C_c^\infty(\Ri^d)$ is a core of $H$ this bound extends to all $\varphi\in D(h)$.
Hence $\one_\Omega D(h)\subseteq D(H^{1/2})=D( h)$.
Then Condition~\ref{tsep2-1} follows from \cite{FOT}, Theorem~1.6.1.
\hfill$\Box$

\bigskip

Theorem~\ref{tsep1} is now a corollary of Theorem~\ref{tsep2} and the following proposition.
\begin{prop}\label{pinv2}
Assume $c_{ij}\in W^{1,\infty}(\Ri^d)$.
Let $\Omega$ be a measurable subset of $\Ri^d$.

Then the following conditions are equivalent.
\begin{tabel}
\item\label{pinv2-1}
$T^{(i)}_tL_2(\Omega)=L_2(\Omega)$ for all $i\in\{1,\ldots,d\}$ and $t\in\Ri$. 
\item\label{pinv2-2}
$T^\psi_tL_2(\Omega)=L_2(\Omega)$ for all $\psi\in C_c^\infty(\Ri^d)$ and $t\in\Ri$. 
\end{tabel}
\end{prop}
\proof\
Note that we have omitted the index used in Section~\ref{S2} to denote  the  space on which the flows act but this should not cause confusion as the following argument only involves $L_2(\Ri^d)$ and $L_\infty(\Ri^d)$.

First,  Condition~\ref{pinv2-1} is equivalent to the conditions $T^{(i)}_t\one_\Omega=\one_\Omega$  for all $i\in\{1,\ldots,d\}$ and $t\in\Ri$ by Lemma~\ref{linv1}.
But the latter conditions are equivalent to 
\begin{equation}
(Y_i^*\varphi,\one_\Omega)=0
\label{elc3.41}
\end{equation}
for all $\varphi\in C_c^\infty(\Ri^d)$ and  $i\in\{1,\ldots,d\}$ by Corollary~\ref{cinv}.

Secondly,  Condition~\ref{pinv2-2}  is equivalent to
\begin{equation}
(Y_\psi^*\varphi,\one_\Omega)=0
\label{elc3.42}
\end{equation}
for all $\varphi, \psi \in C_c^\infty(\Ri^d)$ by similar reasoning.

Now if $\psi\in C_c^\infty(\Ri^d)$ then
\[
(Y_\psi^*\varphi,\one_\Omega)=\sum_{j=1}^d(Y^*_j(\varphi\,\partial_j\psi),\one_\Omega)
\]
for all $\varphi\in C_c^\infty(\Ri^d)$.
Hence (\ref{elc3.41}) implies (\ref{elc3.42}).
Conversely if for $\varphi\in C_c^\infty(\Ri^d)$ and $i\in\{1,\ldots,d\}$ one chooses $\psi\in C_c^\infty(\Ri^d)$  such that $\psi(x)=x_i$ for $x\in \supp\varphi$ then
 (\ref{elc3.42}) implies~(\ref{elc3.41}).
\hfill$\Box$

\section*{Acknowledgements}
The authors would like to thank Ola Bratteli for many discussions on related material.

\end{document}